\documentclass[a4paper,12pt]{amsart}

\usepackage{amssymb,amsfonts,amsmath,amsthm}
\usepackage{amsrefs,mathtools,stmaryrd}
\usepackage{xcolor}
\usepackage{tikz}
\usepackage{verbatim}
\usepackage{doi}
\usepackage{etoolbox}

\usepackage[margin=1in]{geometry}

\newcommand{\Z}{{\mathbb{Z}}}

\newcommand{\R}{{\mathbb{R}}}

\newcommand{\eps}{\varepsilon}
\renewcommand{\phi}{\varphi}

\theoremstyle{plain}
\newtheorem{theorem}{Theorem}

\title{Sharp Gaussian decay for the one-dimensional harmonic oscillator}

\author{Danylo Radchenko} 
\author{Jo\~ao P. G. Ramos}

\begin{document}
\begin{abstract}
We prove a conjecture by Vemuri \cite{Vemuri} by proving  sharp bounds on $\ell^{\kappa}$ sums of Hermite functions multiplied by an exponentially decaying factor. More explicitly, we prove that, for each $y>0,$ we have 
\[
\sum_{n \ge 1} |h_n(x)|^{\kappa} \frac{e^{-\kappa n y}}{n^{\beta}} \ll_y x^{\frac{1}{2} - 2\beta} e^{-\kappa x^2 \tanh(y)/2}, \, \forall \, x \in \R \text{ sufficiently large.}
\]
Our proof involves the classical Plancherel-Rotach asymptotic formula for Hermite polynomials and a careful local analysis near the maximum point of such a bound. 
\end{abstract}

\maketitle

\section{Introduction}

Let $f \in L^2(\R).$ The classical \emph{Hardy uncertainty principle} states that, if 
\begin{equation}\label{eq:Hardy-decay} 
|f(x)| \le C e^{-a \pi |x|^2} \text{ and } |\widehat{f}(\xi)| \le C e^{-b\pi|\xi|^2},
\end{equation} 
for some $C>0,$ almost every $x, \xi \in \R,$ and $a \cdot b > 1,$ then $f \equiv 0,$ where $\widehat{f}$ denotes the Fourier transform of $f.$ Moreover, if $a \cdot b = 1,$ then $f(x) = c \cdot e^{-a \pi |x|^2}$ almost everywhere. 

In the \emph{subcritical} case, that is, when $a \cdot b < 1,$ one has an infinite-dimensional vector space of functions satisfying these conditions, but an important question nonetheless is whether further fine properties of such functions can be derived. This was the main subject of the work \cite{Vemuri} by M. Vemuri, where the author proves that, if a function $f$ satisfies \eqref{eq:Hardy-decay} for $a=b <1,$ then the coefficients in its Hermite expansion must satisfy an \emph{exponential decay}.

More concretely, let $h_n(x)$ be the $n$-th Hermite function, defined by 
	\[h_n(x)=\frac{e^{-x^2/2}H_n(x)}{\sqrt{2^{n}\pi^{1/2}n!}}\,,\]
where $H_n(x)=(-1)^ne^{x^2}\frac{d^n}{dx^n}e^{-x^2}$ is the $n$-th Hermite polynomial. Then any $f \in L^2$ may be written as 
\[
f(x) = \sum_{n \ge 0} c_n h_n(\sqrt{2\pi} x),
\]
where $\|f\|_2^2 = \sum_{n \ge 0} |c_n|^2,$ and the sum converges in the $L^2$ sense. If, moreover, $f$ satisfies \eqref{eq:Hardy-decay} for $a=b < 1,$ then Vemuri showed in \cite[Theorem~2.2]{Vemuri} that
\[
|c_n| \le \frac{C e^{-\alpha n}}{n^{1/4}}, 
\]
where $a = \tanh(2\alpha).$ The author then uses this to conclude a decay estimate for harmonic oscillators: if $\Phi_f$ satisfies 
\begin{equation}\label{eq:harmonic-oscillator}  
i \partial_t \Phi_f(x,t) = \left( \frac{\partial^2}{\partial x^2} - 4\pi^2 x^2 \right) \Phi_f(x,t),
 \end{equation}
with $\Phi_f(x,0) = f(x),$ then, if $f$ satisfies \eqref{eq:Hardy-decay} with $\tanh(2\alpha)=a=b<1,$ then, for any $\beta < \alpha$ and any $t \in \R,$
\[
|\Phi_f(x,t)| \le C_\beta e^{-\tanh(\beta)\pi |x|^2}.
\]
Vemuri, however, conjectured still that an improvement over such a bound must hold: one should be able to replace $\tanh(\beta)$ by $\tanh(\alpha).$ 

To our best knowledge, this conjecture remained almost untouched for many years, since recently A. Kulikov, L. Oliveira and the second author \cite{Kulikov-Oliveira-Ramos} proved Vemuri's conjecture \emph{except} for times $t \in \left\{ \frac{2k+1}{16}, k \in \Z \right\}.$ We refer the reader to \cite{Ga-Ta} and the references therein for further interesting directions. 

The main purpose of this note is to prove Vemuri's conjecture in the remaining cases. That is, 

\begin{theorem}\label{thm:main} Let $E^{\infty}_a = \{ f \in L^2 \colon |f(x)| \ll e^{-a\pi|x|^2}, |\widehat{f}(\xi)| \ll e^{-a\pi|\xi|^2}\}.$ Let $\Phi_f$ be defined as in \eqref{eq:harmonic-oscillator}. Then, if $f \in E^{\infty}_{\tanh(2\alpha)},$ one has
\[
|\Phi_f(x,t)| \le C e^{-\tanh(\alpha)\pi|x|^2}, 
\]
for all $x\in \R$ and $t \in \R.$ 
\end{theorem} 

The main tool used in the proof of Theorem \ref{thm:main} is a decay estimate on sums of Hermite functions with exponential factors, as highlighted by the following result. 

\begin{theorem}\label{thm:decay-hermite}
Let $\kappa > 0$ and $\beta \in \R.$ Then, for any $y>0,$ one has
\begin{equation}\label{eq:bound-hermite-main}  
\sum_{n\ge 1}\frac{|h_n(x)|^{\kappa}}{n^{\beta}}e^{-\kappa ny} \ll_{y,\kappa,\beta} x^{\frac{1}{2} - 2\beta} e^{-\kappa x^2\tanh(y)/2}\,,\qquad \forall \, x \in \R\setminus [-1,1].
\end{equation} 
Moreover, the above bound is sharp, in the sense that the reverse bound also holds whenever $x$ is sufficiently large. Finally, if $\beta \ge \frac{1}{4},$ then the bound above can be extended to \emph{all} $x \in \R.$
\end{theorem}

The proof of Theorem \ref{thm:decay-hermite} involves properly identifying when one can simply use an uniform bound on Hermite functions, and when one needs a more precise asymptotic formula. A similar idea to our proof, used in a different yet related context, can be found in \cite{Boyd1} (see also \cite{Boyd2} for related work). In those works, since the author does not need precise estimates, the usage of an asymptotic formula can be foregone. In our case, however, since we need more precise asymptotics, the Plancherel-Rotach asymptotic formula (see \cite{Sz}) suits our purposes perfectly, and it also shows the sharpness of Theorem \ref{thm:decay-hermite} as we shall see below. 

\subsection*{Notation} We shall denote by $C,c$ absolute constants -- which may change from line to line -- whose exact values are not important for our proof. We also use the notation $a \ll b$ to denote the existence of an absolute constant $C>0$ such that $a \le C b;$ whenever the constant $C$ in such an inequality depends on a parameter $\tau,$ we shall write $a \ll_{\tau} b.$ We shall furthermore write $a = b + O_\tau(c)$ if $|a-b| \ll_{\tau} c.$

\section{Proof of Theorems \ref{thm:main} and \ref{thm:decay-hermite}} 

\subsection{Proof of Theorem \ref{thm:decay-hermite}} We begin with the proof of the second result listed above, since that result will be pivotal for our proof of Theorem \ref{thm:main}. 

We will assume that $x$ is sufficiently large (possibly depending on $y$), and set $N=\max(\lfloor x^2\frac{\tanh(y)}{2y}\rfloor-1,1)$. 
Then, using the fact that $|h_n(x)|\le \pi^{-1/4}$, we get
	\[\sum_{n\ge N}\frac{|h_n(x)|^{\kappa}}{n^{\beta}}e^{-\kappa ny} \ll_y x^{-2\beta} e^{-\kappa x^2\tanh(y)/2} \,,\]
so that we only need to deal with the sum for $1\le n\le N-1$. Note that 
$N+1\le (1-\eps)\frac{x^2}{2}$ for some $\eps>0$ depending on $y$ 
(e.g. take $\eps = \frac{1}{2}(1-\tanh(y)/y)$), so that we are in the monotonic regime for 
Hermite functions. In this case the Plancherel-Rotach asymptotic formula (see~\cite{Pl-Ro}, and also~\cite[(8.22.13)]{Sz}) provides an upper bound
	\[|h_n(x)| \ll \frac{\exp((n+\frac12)\phi_n-(n+1)\sinh(\phi_n)\cosh(\phi_n))}{2^{3/4}\pi^{1/2}n^{1/4}\sinh(\phi_n)^{1/2}}
	\qquad 2\le n+1 \le (1-\eps)\frac{x^2}{2}\,,\]
where $\phi_n$ is defined by $x=\sqrt{2(n+1)}\cosh(\phi_n)$; that is, $\varphi_n = \cosh^{-1}(x/\sqrt{2(n+1)}).$ Thus, for $1\le n<N$, using the fact that in this range $\cosh(\phi_n)\ge (1-\eps)^{-1/2}$ and the relation $\cosh^2 - \sinh^2 = 1$ repeatedly, we get
\begin{equation}\label{eq:bound-hermite} 
\frac{|h_n(x)|^{\kappa} e^{-\kappa ny}}{n^{\beta}}\ll_{y,\kappa} n^{-1/4 - \beta}\exp\Big(\kappa \left(n\phi_n-ny-\frac{x}{2}\sqrt{x^2-2(n+1)}\right)\Big)\,.
\end{equation}

We then define the argument function for the exponential in the bound above as
$$A(n) = n \varphi_n - ny - \frac{x}{2} \sqrt{x^2 - 2(n+1)}.$$ 
It follows then that $A$ is a smooth function of $n$ (now regarded as a real variable in the interval $[1,N-1]$), and hence a quick computation shows that
	\begin{align*}
	\partial_n A(n)  &= \phi_n+\frac{x}{2(n+1)\sqrt{x^2-2n-2}}-y\,,\\
	\partial_n^2 A(n) &= x\frac{(2n+5)(n+1)-x^2(n+2)}{2(n+1)^2(x^2-2n-2)^{3/2}}\,.
	\end{align*}
We now claim that $\partial_n^2 A(n)< -\frac{c}{x^2}$ for $1<n<\frac{x^2-4}{2}$, for some absolute constant $c>0$. As a matter of fact, if we call $t_n = x^2/(2n+2),$ this assertion is equivalent to 
\begin{equation}\label{eq:bound-second-derivative} 
\left( 1 + \frac{1}{n+1} \right) t_n - \left(1+ \frac{3}{2n+2}\right) \ge c \left( 1 - \frac{1}{t_n}\right)^{3/2},
\end{equation}
for all $n \in (1,\frac{x^2-4}{2})$ and some constant $c>0.$ One readily sees that, for $c$ sufficiently small, the derivative of the difference between the two sides of \eqref{eq:bound-second-derivative} with respect to $t_n$ is nonnegative, and since the inequality holds if $t_n = 1 + \frac{1}{n+1},$ the assertion is proved. 

This implies, on the other hand, that $A(n)$ has a unique maximum on $[1,N-1]$. We claim that for all sufficiently large $x$ the point $n_{\max}$ where $A(n)$ is maximized satisfies
	\[n_{\max}=\frac{x^2}{2\cosh^2(y)}+O_y(1)\,.\]
Indeed, the condition $\partial_n A(n)=0$ can be rewritten as
	\[\phi_n+\frac{\cosh(\phi_n)^3}{x^2\sinh(\phi_n)} = y\,.\]
Let us look then at the solutions $u_x$ to the equation 
\[
u_x + \frac{\cosh(u_x)^3}{x^2 \sinh(u_x)} = y.
\]
Notice first that $u_x > 0$ since $y>0.$ This implies that $u_x \le y,$ and thus $u_x  =y + O_y(x^{-2}).$ Using this fact with $\varphi_n$ in place of $u_x,$ together with the definition of $\varphi_n,$ yields the claim on the asymptotic form of $n_{\max}.$ Plugging $n_\text{max}$ into $A$ then yields
	\begin{align*}
	\max_{1\le n\le N-1}A(n) 
	&\le A(n_{\max})=n_{\max}(\varphi_{n_{\max}} - y) -(n_{\max}+1)\sinh(\phi_n)\cosh(\phi_n)\\
	&= -\frac{x^2\sinh(y)\cosh(y)}{2\cosh^2(y)} + O_y(1)
	 = -\frac{x^2\tanh(y)}{2} + O_y(1)\,.
	\end{align*}
Finally, using the fact that $\partial_n^2 A(n)<-c/x^2$, we get
	\[A(n)\le A(n_{\max})-\frac{c(n-n_{\max})^2}{2x^2}\,,\]
and so, if we define $\lambda = \frac{n_{\max}}{x^2},$ we get
	\begin{align*}
	\sum_{n=1}^{N-1}\frac{\exp(\kappa A(n))}{n^{1/4 +\beta}} &\ll_{y,\kappa} e^{-\kappa x^2\tanh(y)/2}\sum_{n=1}^{N-1}\frac{\exp(-\kappa c(n-\lambda x^2)^2/x^2)}{n^{1/4 + \beta}}\\
	&\ll_{y,\kappa,\beta} e^{-\kappa x^2\tanh(y)/2}\Big(\sum_{n<\lambda x^2/2}
	\frac{e^{-\kappa c\lambda^2 x^2/4}}{n^{1/4 + \beta}}+\sum_{N>n>\lambda x^2/2}\frac{e^{-\kappa c(n-\lambda x^2)^2/x^2}}{n^{1/4 + \beta}}\Big)\\
	&\ll_{y,\kappa,\beta} e^{-\kappa x^2\tanh(y)/2}\Big(\max\{1,x^{\frac{3}{2} - 2\beta}\}e^{-c\lambda^2 x^2/4}+\frac{\sum_{n \in \Z}e^{-c'n^2/x^2}}{x^{\frac{1}{2} + 2\beta}}\Big). \qedhere
	\end{align*}
Now, we note that by Poisson summation the function $g_\delta(x) := \sum_{n \in  \Z} e^{-\delta \pi n^2/x^2}$ satisfies 
    \[
    g_\delta(x) = \frac{x}{\sqrt{\delta}} \sum_{n \in \Z} e^{-\pi n^2 x^2/\delta}.
    \]
Since the series on the right-hand side above is uniformly bounded for $x$ sufficiently large, we get that $g(x)/x$ is bounded for $x \in \R.$ Applying this to $\delta = \frac{c'}{\pi}$ finishes the proof of \eqref{eq:bound-hermite-main}.

\vspace{3mm}
We shall now see that Theorem \ref{thm:decay-hermite} is sharp by proving that the inequality can be reversed for $x$ sufficiently large. 

As a matter of fact, the Plancherel-Rotach asymptotic formula, applied in the form of \cite[(8.22.13)]{Sz}, implies that the inequality in \eqref{eq:bound-hermite} can be \emph{reversed}, as long as $\frac{2n+2}{x^2}$ belongs to an interval $I$ strictly contained in $(0,1).$ It then follows that, in that case, also \eqref{eq:bound-second-derivative} can be reversed, as long as one takes $c$ to depend on the interval $I.$ This implies that $\partial_n^2 A(n) \ge -\frac{c'}{x^2}$ for such $n,$ and hence also $A(n) \ge A(n_{\max}) - \frac{c'(n-n_{\max})^2}{2x^2}$ in that range. Therefore,
\begin{align*}
\sum_{n \colon \frac{2n+2}{x^2} \in I} |h_n(x)|^{\kappa} \frac{e^{-\kappa n y}}{n^{\beta}} & \gg_{y,\kappa,\beta,I} \sum_{n \colon \frac{2n+2}{x^2} \in I} \frac{\exp(\kappa A(n))}{n^{1/4 +\beta}} \cr 
    & \gg_{y,\kappa,\beta,I} e^{-\kappa x^2\tanh(y)/2}\Big(\sum_{n \colon \frac{2n+2}{x^2} \in I}\frac{e^{-\kappa c'(n-\lambda x^2)^2/x^2}}{n^{1/4 + \beta}}\Big) \cr 
    & \gg_{y,\kappa,\beta,I} x^{-\frac{1}{2} - 2\beta} e^{-\kappa x^2\tanh(y)/2} \Big( \sum_{n \colon \frac{2n+2}{x^2} \in I} e^{-\kappa c'(n-\lambda x^2)^2/x^2} \Big).
\end{align*}
Choose then $I = (\lambda/2,\tanh(y)/y).$ Thus, one has
\begin{align*} 
\sum_{n \colon \frac{2n+2}{x^2} \in I} e^{-\kappa c'(n-\lambda x^2)^2/x^2} & \ge \sum_{m=0}^{\lfloor \left(\frac{\tanh(y)}{y} - \lambda\right)x^2 \rfloor} e^{-\kappa c'' m^2/x^2} \gg \int_0^{\frac{1}{2} \left(\frac{\tanh(y)}{y} - \lambda\right)x^2 } e^{-\kappa c'' t^2/x^2} \, dt \gg_{\kappa,y} x. 
\end{align*}
Putting these together, one obtains that 
\[
\sum_{n \ge 0} |h_n(x)|^{\kappa} \frac{e^{-\kappa n y}}{n^{\beta}} \gg_{y,\kappa,\beta} x^{\frac{1}{2} - 2\beta} e^{-\kappa x^2 \tanh(y)/2},
\]
which concludes the proof of the claim, and hence also the proof of Theorem \ref{thm:decay-hermite}. 

\subsection{Proof of Theorem \ref{thm:main}} We now use Theorem \ref{thm:decay-hermite} to prove Theorem \ref{thm:main}. Note that the function $\Phi_f(x,t)$ can be \emph{explicitly} written in terms of the Hermite basis; indeed, if
\[
f(x) = \sum_{n \ge 0} c_n h_n(\sqrt{2\pi}x),
\]
then 
\[
\Phi_f(x,t) = \sum_{n \ge 0} e^{2(2n+1)\pi i t} c_n h_n(\sqrt{2\pi}x). 
\]
Moreover, we recall once more that Vemuri proved in \cite[Theorem~2.2]{Vemuri} that, if $f \in E^{\infty}_{\tanh(2\alpha)},$ then $|c_n| \ll \frac{e^{-\alpha n}}{n^{1/4}}.$ Thus, 
\[
|\Phi_f(x,t)| \ll \sum_{n \ge 0} |h_n(\sqrt{2\pi}x)| \frac{e^{-\alpha n}}{n^{1/4}}.
\]
Applying Theorem \ref{thm:decay-hermite} with $\kappa =1, \beta = 1/4$ and $y = \alpha$ then implies that 
\[
|\Phi_f(x,t)| \ll_{\alpha} e^{-\tanh(\alpha)\pi |x|^2},
\]
as stated in \cite[Conjecture~3.2]{Vemuri}, which concludes the proof of Theorem \ref{thm:main}. 

\section*{Acknowledgements} 

J.P.G.R. acknowledges financial support by the European Research Council under the Grant Agreement No. 721675 “Regularity and Stability in Partial Differential Equations (RSPDE)”. D.R. acknowledges financial support by the European Research Council (ERC Starting Grant No. 101078782).

\end{document}